\documentclass[11pt,reqno]{amsart}
\usepackage{amssymb, amsmath, amsthm,amsrefs,graphicx,marvosym}
\usepackage{pstricks}
\usepackage{latexsym}
\usepackage{color}
\usepackage{exscale}
\usepackage{bm, bbm, mathrsfs}
\usepackage{hyperref}

\newtheorem*{theorem*}{Theorem}
\newtheorem*{lemma*}{Lemma}

\newtheorem{theorem}{Theorem}
\newtheorem{lemma}{Lemma}[section]

\newtheorem{definition}[lemma]{Definition}
\newtheorem{remark}[lemma]{Remark}
\newtheorem{example}[lemma]{Example}

\DeclareMathOperator{\spa}{span}

\newcommand{\RR}{\mathbb R}

\def\s0{{ s_0}}

\def\ts0{{\tilde s_0}}

\def\eq#1{(\ref{#1})}
\def\nn{\nonumber}

\def\({\left(\begin{array}{cccccc}}
\def\){\end{array}\right)}

\def\bes{\begin{eqnarray}}
\def\ees{\end{eqnarray}}

\newcommand{\del}{\partial}

\newcommand{\beq}{\begin{equation}}
\newcommand{\eeq}{\end{equation}}
\newcommand{\bea}{\begin{eqnarray}}
\newcommand{\eea}{\end{eqnarray}}
\newcommand{\baln}{\begin{align}}
\newcommand{\ealn}{\end{align}}
\newcommand{\beann}{\begin{eqnarray*}}
\newcommand{\eeann}{\end{eqnarray*}}

\newcommand{\R}{\ensuremath{\mathfrak{R}}}

\newcommand{\Rn}{\mathbb R^n}
\newcommand{\Rm}{\mathbb R^m}
\newcommand{\og}{\bar g}
\newcommand{\rset}{\{\br_i\}}
\newcommand{\Xia}{\Xi_\alpha}
\newcommand{\ga}{g_\alpha}
\newcommand{\ua}{u_\alpha}
\newcommand{\atx}{\big|_x}
\newcommand{\fa}{f^\alpha}
\newcommand{\fia}{f_i^\alpha}
\newcommand{\Ia}{I_\alpha}
\newcommand{\ia}{i_\alpha}
\newcommand{\pa}{{p_\alpha}}
\newcommand{\wtO}{\widetilde \Omega}
\newcommand{\wtU}{\widetilde \Upsilon}
\newcommand{\mC}{\mathcal C}
\newcommand{\psia}{\psi_\alpha}
\newcommand{\xia}{\xi_\alpha}
\newcommand{\wua}{{\widetilde u}_\alpha}
\newcommand{\U}{\mathcal U}
\newcommand{\ijj}{{i_j}}
\newcommand{\ikk}{{i_k}}

\newcommand{\br}{\mathbf r}

\newcommand{\wu}{{\tilde u}}

\newcommand{\ox}{{\bar x}}
\newcommand{\dthm}{Th\'eor\`eme }
\newcommand{\dthms}{Th\'eor\`emes }

\newcommand{\pf}{\begin{proof}}
\newcommand{\foorp}{\end{proof}}

\newcommand{\Addresses}{{
  \bigskip
  \footnotesize

  M.~Benfield, \textsc{La Mesa, CA 
} \texttt{(mike.benfield@gmail.com).}

  \medskip

 H.~K.~Jenssen, \textsc{Department of
Mathematics, Pennsylvania State University}
 \texttt{(jenssen@math.psu.edu).}

  \medskip

  I.~A.~Kogan, \textsc{Department of Mathematics, North Carolina State University}\texttt{(iakogan@ncsu.edu).}

}}

\begin{document}

\title{A generalization of an integrability theorem of Darboux}
\author{Michael Benfield}

\author{Helge Kristian Jenssen }

\author{Irina A.\ Kogan}

\date{\today}

\begin{abstract}
In his monograph ``Syst\`emes Orthogonaux'' \cite{dar} Darboux
stated three theorems providing local existence 
and uniqueness of solutions to first order systems of the type  
\[\partial_{x_i} u_\alpha(x)=f^\alpha_i(x,u(x)),
\quad i\in I_\alpha\subseteq\{1,\dots,n\}.\] 
For a given point $\bar x\in \RR^n$ it is assumed that the values of 
the unknown $u_\alpha$ are given locally near $\bar x$ along 
$\{x\,|\, x_i=\bar x_i \, \text{for each}\, i\in I_\alpha\}$.
The more general of the theorems, \dthm III, was proved by Darboux
only for the cases $n=2$ and $3$.

In this work we formulate and prove a generalization of Darboux's
\dthm III which applies to systems of the form 
\[{\mathbf r}_i(u_\alpha)\big|_x = f_i^\alpha (x, u(x)), \quad i\in
I_\alpha\subseteq\{1,\dots,n\}\] where $\mathcal R=\{{\mathbf r}_i\}_{i=1}^n$ 
is a fixed local frame of vector fields near $\bar x$. The data for $u_\alpha$
are prescribed along a manifold $\Xi_\alpha$ containing $\bar x$ and 
transverse to the vector fields $\{{\mathbf r}_i\,|\, i\in I_\alpha\}$. 
We identify a certain Stable Configuration Condition (SCC). This is a geometric condition 
that depends on both 
the frame $\mathcal R$ and on the manifolds $\Xi_\alpha$; it is
automatically met in the case considered by Darboux \cite{dar}.
Assuming the SCC and the relevant integrability 
conditions are satisfied, we establish local existence and uniqueness 
of a $C^1$-solution via Picard iteration for any number of independent variables $n$. 
\end{abstract}

\maketitle


\tableofcontents

\section{Introduction}\label{intro}
Darboux, in Chapitre I, Livre III in his monograph ``Syst\`emes Orthogonaux''
 \cite{dar}, stated three integrability theorems (``\dthms I-III'') for certain types of first order
systems of PDEs.
The theorems apply to systems of the form
\begin{equation}\label{stand_eqns}
	\partial_{ x_i}u_\alpha(x)=f^\alpha_i(x,u(x))
\end{equation}
where $u=(u_1, \ldots, u_m)$ denotes the vector of unknown functions, the
independent variables $x=(x_1, \ldots, x_n)$ range over an open set about a
fixed point $\bar x \in \mathbb R^n$, and the $f^\alpha_i$ are given
$C^1$-functions from an appropriate open subset of $\RR^{n+m}$ to $\RR$. The data for these systems consist of $C^1$-functions $g_\alpha$
prescribing the unknowns $u_\alpha$ along certain affine subspaces
through the point $\bar x \in \mathbb R^n$.

For each $\alpha=1,\dots,m$, we let $I_\alpha$ be the set of all indices 
$i\in\{1,\dots,n\}$ for which the system contains the equation $\partial_{ x_i}u_\alpha=f^\alpha_i(x,u(x))$.
The cases covered by Darboux's three theorems can then be described 
as follows:

\dthm I applies to determined systems, which are characterized by the 
        requirement that $|I_\alpha|=1$ for all $\alpha=1,\dots,m$. For each 
        $\alpha$, letting $I_\alpha=\{i_\alpha\}$, the data for $u_\alpha$ are 
        prescribed near $\bar x$ along the hyperplane $\{x\,|\,x_{i_\alpha}=\bar x_{i_\alpha}\}$. 
        In the special case that $i_\alpha$ is the same index for all $\alpha$, Darboux's 
        \dthm~I reduces to the standard local existence and uniqueness result for ODEs
        with parameters (see \cite{hart}).
        
\dthm II is the PDE version of Frobenius' theorem for completely 
        overdetermined systems. This situation is characterized by $|I_\alpha|=n$ 
        for all $\alpha=1,\dots,m$, i.e.\  the derivatives of all unknowns are prescribed 
        in all coordinate directions. In this case, the data prescribes the value of each $u_\alpha$ 
        at the point $\bar x$. The integrability conditions require 
        the partial derivatives given by the system to be consistent with equality of 2nd
        order mixed partial derivatives. Under these conditions Darboux's \dthm II guarantees
        a unique local solution of the PDE system with the assigned data.
        
\dthm III applies to the general case where the elements, as well 
        as the cardinality, of $I_\alpha$ may vary with $\alpha$. 
        The data are assigned as follows: if 
	$I_\alpha=\{i_1,\dots,i_{p_\alpha}\}$, we prescribe a function $g_\alpha$ 
	along the affine subspace 
	$\Xi_\alpha:=\{x\,|\, x_{i_j}=\bar x_{i_j},\ 1\leq j\leq p_\alpha\}$, and 
	require that $u_\alpha|_{\Xi_\alpha}=g_\alpha$.
	Under the appropriate integrability conditions, detailed below in Section~\ref{setup}, 
	Darboux's \dthm~III guarantees
        a unique local solution of the PDE system with the assigned data.

We note that \dthms I and II are particular cases of \dthm III.
Another special case of systems, for which the index sets $I_\alpha$ are the same for 
all $\alpha$, was addressed separately in \cite{mike_thesis}.

Darboux stated his \dthm III for any number of independent 
variables. He provided a proof only in the cases 
with two and three independent variables ($n=2$ or $3$), which sufficed for 
his investigation of triply orthogonal systems in \cite{dar}.
For $n=2$ his proof of \dthm III used \dthm I; for $n=3$ he used both the result for $n=2$
as well as his \dthm I. 
Given Darboux's partial proof, it is natural to try and establish his \dthm III via  
induction on the number  $n$ of independent variables. While this is possible, we 
have been able to do so only through an involved, combinatorial  
argument (see our unpublished note \cite{bjk_unpub}). Furthermore, 
this inductive approach does not apply to the more general situation 
we consider in the present paper.
Instead, we shall provide a direct proof that applies to  more general systems 
with any number of independent variables.  


 \noindent Our results generalize Darboux's \dthm III in two ways:
\begin{enumerate}
 \item[(i)]  The unknowns may be differentiated along vector fields in 
	a  fixed frame $\mathcal R=\{\br_i\}_{i=1}^n$ defined 
	near $\bar x$. That is, for each $\alpha=1,\dots,m$, there is an 
	index set $I_\alpha\subseteq \{1,\dots,n\}$ such that the system 
	contains the equations 
	\beq\label{gen_eqns}
		\br_i(u_\alpha)\big|_x = f_i^\alpha (x, u(x))\qquad\text{for each $i\in I_\alpha$.}
	\eeq
	As in Darboux's \dthm III, the elements and cardinality of the index sets 
	$I_\alpha$ may vary with $\alpha$. 
    
  \item[(ii)]  The prescribed data $g_\alpha$ for the unknown $u_\alpha$ may be given 
	along a manifold $\Xi_\alpha$ through the point $\bar x$ which is transverse to the 
	vector fields $\br_i$ with $i\in I_\alpha$.
\end{enumerate}

The claim is that, under the appropriate integrability conditions
(generalizing those of Darboux's \dthm III), the PDE system \eq{gen_eqns} has
a unique local solution which takes on the assigned data. A precise
formulation is provided in our Theorem~\ref{thm}  in Section~\ref{sec:theorem}.

However, our proof requires what we refer to as a 
{\em Stable Configuration Condition} (SCC) to be satisfied. 
The formulation of the SCC is somewhat technical (see Section 
\ref{scc_cond} and Definitions \ref{def-access} and \ref{def-stable} 
below). Roughly speaking, this condition is required to guarantee 
that the natural Picard iteration scheme is well defined. We note 
that the validity of the SCC depends on both the frame $\mathcal R$ 
and on the relative location of the manifolds $\Xi_\alpha$ that carry 
the data; see Section \ref{scc_cond} below for a concrete example.
Also, it is immediate to verify that the SCC is met in the setting of 
Darboux's original treatment where 
$\br_i\equiv \del_{x_i}$, $i=1,\dots,n$, and $\Xi_\alpha=\{x\,|\, x_i=\bar x_i, 
\, i\in I_\alpha\}$.

Concerning regularity, we assume that the frame $\{\br_i\}_{i=1}^n$, the
functions $f^\alpha_i$, the manifolds $\Xi_\alpha$, and the data 
$g_\alpha$ are all $C^1$-smooth. A
solution refers to a $C^1$-smooth function $u=(u_1,\dots,u_m)$ which satisfies
the PDEs and the data in a classic, pointwise manner on a neighborhood of the
given point $\bar x$.

The rest of the present paper is organized as follows. In Section \ref{setup}
we review Darboux's original \dthm III and the partial proof provided by
Darboux. We also indicate how our approach in this paper differs from that of
Darboux. Section \ref{scc_cond} considers a simple system of equations to highlight the 
role of the Stable Configuration Condition (SCC): for a determined system of 
two equations for two unknowns in the plane, we show how the
relative location of the two data manifolds $\Xi_1$ and $\Xi_2$ can
yield radically different behavior in terms of the domains of definition of the natural 
Picard iterates. Finally, in Section \ref{sec:theorem} we formulate and prove our Theorem.
A key part of the proof is a technical lemma about the 
``restricted'' system obtained by considering the same set of equations 
as in the original system, but restricted to certain sub-manifolds defined 
in terms of the frame vector fields $\br_i$; see Lemma~\ref{lemma} below.

\section{Review of Darboux's work and the Stable Configuration Condition (SCC)}
\subsection{Darboux's setup and result}\label{setup}
We first consider the situation addressed by Darboux in his \dthm III: for 
each unknown $u_\alpha$, the system consists of the equations
\beq\label{stand_syst}
	\partial_{x_i} u_\alpha(x)=f^\alpha_i(x,u(x))
	\qquad \text{for $i\in I_\alpha\subseteq \{1,\dots,n\}$},
\eeq
where $f^i_\alpha$ are $C^1$- smooth functions on $\RR^{n+m}$.
Setting
\beq\label{data_mf}
	\Xi_\alpha:=\{x\,|\, x_i=\bar x_i \text{ for } i\in I_\alpha\},
\eeq
we prescribe the data
\beq\label{stand_data}
	u_\alpha\big|_{\Xi_\alpha}=g_\alpha,
\eeq
for a given $C^1$- smooth function $g_\alpha:\Xi_\alpha\to \RR$.

Next consider the integrability conditions which need to be imposed.
Let $u_\alpha$ be an unknown for which the system prescribes two 
distinct partial derivatives, say
\[\partial_{x_i} u_\alpha(x)=f^\alpha_i(x,u(x))
\qquad \text{and}\qquad
\partial_{x_j} u_\alpha(x)=f^\alpha_j(x,u(x)) \]
where $i\neq j$ and $i,j\in I_\alpha$.
The derivatives prescribed by the system need to be consistent with 
equality of mixed partial derivatives. That is, the expressions
\beq\label{ij_deriv}
	\partial^2_{x_jx_i} u_\alpha(x)
	=\partial_{x_j} f^\alpha_i(x,u(x))
	+\sum_{\beta=1}^m \partial_{u_\beta}f^\alpha_i(x,u(x))\partial_{x_j} u_\beta(x)
\eeq
and 
\beq\label{ji_deriv}
	\partial^2_{x_ix_j} u_\alpha(x)
	=\partial_{x_i} f^\alpha_j(x,u(x))
	+\sum_{\beta=1}^m \partial_{u_\beta}f^\alpha_j(x,u(x))\partial_{x_i} u_\beta(x)
\eeq
should agree. Since the system \eq{stand_eqns} may not prescribe all
the partials $\partial_{x_j} u_\beta$ and $\partial_{x_i} u_\beta$ appearing on 
the right-hand sides of \eq{ij_deriv} and \eq{ji_deriv}, this puts constraints 
on which dependent variables $u_\beta$ the functions $f_i^\alpha$ and 
$f_j^\alpha$ may depend on. This is brought out in the following example
which is the simplest case of an overdetermined system where 
Darboux's \dthm III applies.
\begin{example}\label{3rd_ex1}
	Consider a system of 3 equations for 2 unknowns in 2 independent 
	variables. Let the unknowns be $u$ and $v$, 
	the independent variables be $x$ and $y$, and assume that the 
	equations are
	\begin{align}
		u_x &= f(x,y,u,v)\label{3rd_ex1_1}\\
		v_x &= \phi(x,y,u, v)\label{3rd_ex1_2}\\
		v_y &= \psi(x,y,u,v)\,.\label{3rd_ex1_3}
	\end{align}
	The data in this case take the form
	\begin{align}
		u(\bar x,y) &= g_1(y)\label{3rd_ex1_4}\\
		v(\bar x,\bar y) &= g_2,\label{3rd_ex1_5}
	\end{align}
	where $g_1$ is a given function and $g_2$ is a given constant.
	The integrability condition is imposed to ensure that
         the prescription of the  two  partial 
	derivatives of the unknown $v$  is consistent with the equality 
	of the partial derivatives $(v_x)_y=(v_y)_x$. 
	To derive these conditions we expand $\del_y [\phi(x,y,u,v)]=\del_x [\psi(x,y,u,v)]$
	applying the chain rule, to obtain
	\beq\label{ex-ic}\phi_y+\phi_u \,u_y+\,\phi_v\, v_y=\psi_x+\psi_u\,u_x+\,\psi_v\, v_x.\eeq 
	We next  substitute the derivatives given by the system  \eq{3rd_ex1_1}-\eq{3rd_ex1_3} into \eq{ex-ic}. 
	However,  the system does not provide an expression for $u_y$, and we must therefore impose 
	the condition
	      \beq\label{ex-cond-a}\phi_u=0.\eeq
	All other partial derivatives of $u$ and $v$ appearing in \eq{ex-ic} are prescribed by 
	\eq{3rd_ex1_1}-\eq{3rd_ex1_3}, and we obtain the condition:  
	 \beq\label{ex-cond-b}\phi_y+\,\phi_v\, \psi=\psi_x+\psi_u\,f+\,\psi_v\, \phi .\eeq
	Conditions \eq{ex-cond-a} and \eq{ex-cond-b} comprise the  integrability conditions for the system  \eq{3rd_ex1_1}-\eq{3rd_ex1_3}.
	If these conditions hold as identities in an $(x,y,u,v)$-neighborhood 
	of $(\bar x,\bar y,g_1(\bar y), g_2)$, then
	 Darboux's \dthm~III
	guarantees the existence of a unique local $C^1$-smooth solution $(u(x,y),v(x,y))$ 
	to \eq{3rd_ex1_1}-\eq{3rd_ex1_3} near $(\bar x,\bar y)$ taking on the data 
	\eq{3rd_ex1_4}-\eq{3rd_ex1_5}.\end{example}
In the general setting of the system \eq{stand_syst}, the integrability conditions 
require that, whenever $\alpha\in\{1,\dots,m\}$ and $i,j \in I_\alpha$ with $i\neq j$, 
then the following should hold: for all $\beta\in\{1,\dots,m\}$ with $i\notin I_\beta$ 
we have 
\[\partial_{u_\beta}f_j^\alpha = 0\]
and
\[\del_{x_i} f_j^\alpha +
   { \sum_{\beta: i \in I_\beta} \left(\partial_{u_\beta} f_j^\alpha\right)\, f_i^\beta} 
    \,\equiv\, \del_{x_j} f_i^\alpha 
    + {\sum_{\beta:j\in I_\beta} \left(\partial_{u_\beta} f_i^\alpha\right)\, f_j^\beta}.\]
If these conditions hold as identities in a neighborhood 
of $(\bar x, g(\bar x))$, then Darboux's \dthm~III
guarantees the existence of a unique local $C^1$-smooth solution $u(x)$ 
to \eq{stand_syst} in a neighborhood of  $\bar x$ that takes on the data 
\eq{stand_data}.

Due to the particular structure of the systems under consideration 
(viz., each equation contains a single derivative
for which it is solved), it is natural to base a proof of existence on 
Picard iteration. Indeed, it is immediate to write down a functional 
map for which any solution of \eq{stand_syst}-\eq{stand_data}
must be a fixed point (see \eq{eq:functional} below).

Now, in the particular situation addressed by Darboux's \dthm I (for determined systems), 
one can verify that a fixed point exists and provides a solution 
of the original system. This is how Darboux established his \dthm I. 
However, in the more general situation of overdetermined systems addressed 
by his \dthm III, such an approach appears to be more challenging. This 
circumstance might explain why Darboux \cite{dar} did not provide a general 
proof based directly on Picard iteration for his \dthm III. Instead, 
for $n=2$ case he exploited   \dthm I. 
For  $n=3$ case he   identified sub-systems  that can be treated by   \dthm I or by $n=2$ case. These sub-systems  are  solved  in a ``right'' order so that the solution of one sub-system provides initial data to the next.  Darboux states that the general proof will be too technical and, therefore, he restricts himself to the cases of $n=2$ and $n=3$ as they are sufficient for the applications he considers.\footnote{``Pour \'etablier cette importante proposition, sans employer un trop grand luxe de notations, nous nous bornerons au cas de deux et de trois variables ind\'ependantes, qui  suffira d'ailleurs   pour les applications que nous avons en vue.'' \cite{dar}  p.~336.}  Darboux's treatment of $n=3$ case indicates a possibility of  an inductive proof for an arbitrary $n$, which we accomplish in \cite{bjk_unpub}. The inductive proof turns out, indeed, to be quite  technical.  Moreover,  the same approach can not be applied  
to a more general problem considered here, because even if the initial system satisfies the hypothesis of our Theorem~\ref{thm}, the sub-systems appearing in the inductive proof may not satisfy these hypothesis.

In fact, as we shall show below, it is possible to provide a direct argument
based on Picard iteration also for overdetermined systems of the more general
type \eq{gen_eqns} with any number of independent variables. 
The key observation is that it suffices to consider a
certain ``restricted'' system which consists of the same equations as the
original system, but now required to hold only along certain submanifolds
containing the given point $\bar x$. For this, seemingly weaker, restricted
system, we establish existence of a solution $\tilde u$
via Picard iteration. This result (Lemma~\ref{lemma} below) is the first main new
ingredient in our approach. (We note that this part of the argument does not 
involve any use of the integrability conditions.)

In  the more general setting 
described in  (i) and (ii) above, the argument of Lemma~\ref{lemma} 
becomes complicated by the need to work with different
coordinate systems for different components $u_\alpha$ of the solution $u$.
This is where we have found it necessary to introduce the Stable Configuration
Condition (SCC) illustrated in the next section. We stress that, in the more general setting,  the SCC is relevant 
already for generalizing Darboux's \dthm I  (i.e., the it is not about determinacy
or over-determinacy of the system). However, in the original setting described  by Darboux \cite{dar}, the SCC condition is trivially satisfied.
Therefore, our paper contains a direct proof of  Darboux's \dthm III,  for an arbitrary number of variables.

Before considering an example explaining the relevance of the SCC, 
we outline the last step of the proof: showing that the fixed point $\tilde u$ 
is a solution of the original system on a {\em full} $\RR^n$-neighborhood of $\bar x$.
This is accomplished by showing that the quantities 
\[ A_i^\alpha(x) = \br_i(\tilde u_\alpha)|_x - f_i^\alpha(x, \tilde u(x)),
\qquad 1\leq \alpha\leq m,\,  i \in I_\alpha, \]
satisfy certain linear, homogeneous equations which form a restricted system
of the type covered by our Lemma. Only at this point are the integrability
conditions used. As this latter system admits the trivial solutions
$A_i^\alpha\equiv0$ on a full neighborhood of $\bar x$, it follows from the
uniqueness part of the Lemma that $\br_i(\tilde u_\alpha)|_x = f_i^\alpha(x,
\tilde u(x))$ for all $x$ near $\bar x$, thereby completing the proof.

\subsection{The Stable Configuration Condition (SCC)} \label{scc_cond}
To illustrate the relevance of the SCC we consider the following simple example
of a system for the two unknown scalar functions $u(x,y)$ and $v(x,y)$ 
of two independent variables:
\begin{align}
	u_x&=v\label{u}\\
	v_y&=u.\label{v}
\end{align}
The system \eq{u}-\eq{v} is a determined system of the form \eq{gen_eqns}
with $m=n=2$, $u_1=u$, $u_2=v$, $\br_1=\del_x$, and $\br_2=\del_y$.
We consider the system near the origin in the $(x,y)$-plane and let $M$ and 
$N$ be the two straight lines 
\[M:=\{(x,y)\,|\,x=a y\}\qquad N:=\{(x,y)\,|\,y=b x\},\]
where $0<\frac{1}{a}<b$.
Next, we consider separately the two cases:
\begin{itemize}
	\item[(a)] The data for $u$ are prescribed along $\Xi_1=N$, and the 
	data for $v$ are prescribed along $\Xi_2=M$.
	\item[(b)] Vice versa: $u$ is prescribed along $\Xi_1=M$ and 
	$v$ is prescribed along $\Xi_2=N$.
\end{itemize}
\begin{figure}\label{Figure_1}
	\centering
	\includegraphics[width=11cm,height=8cm]{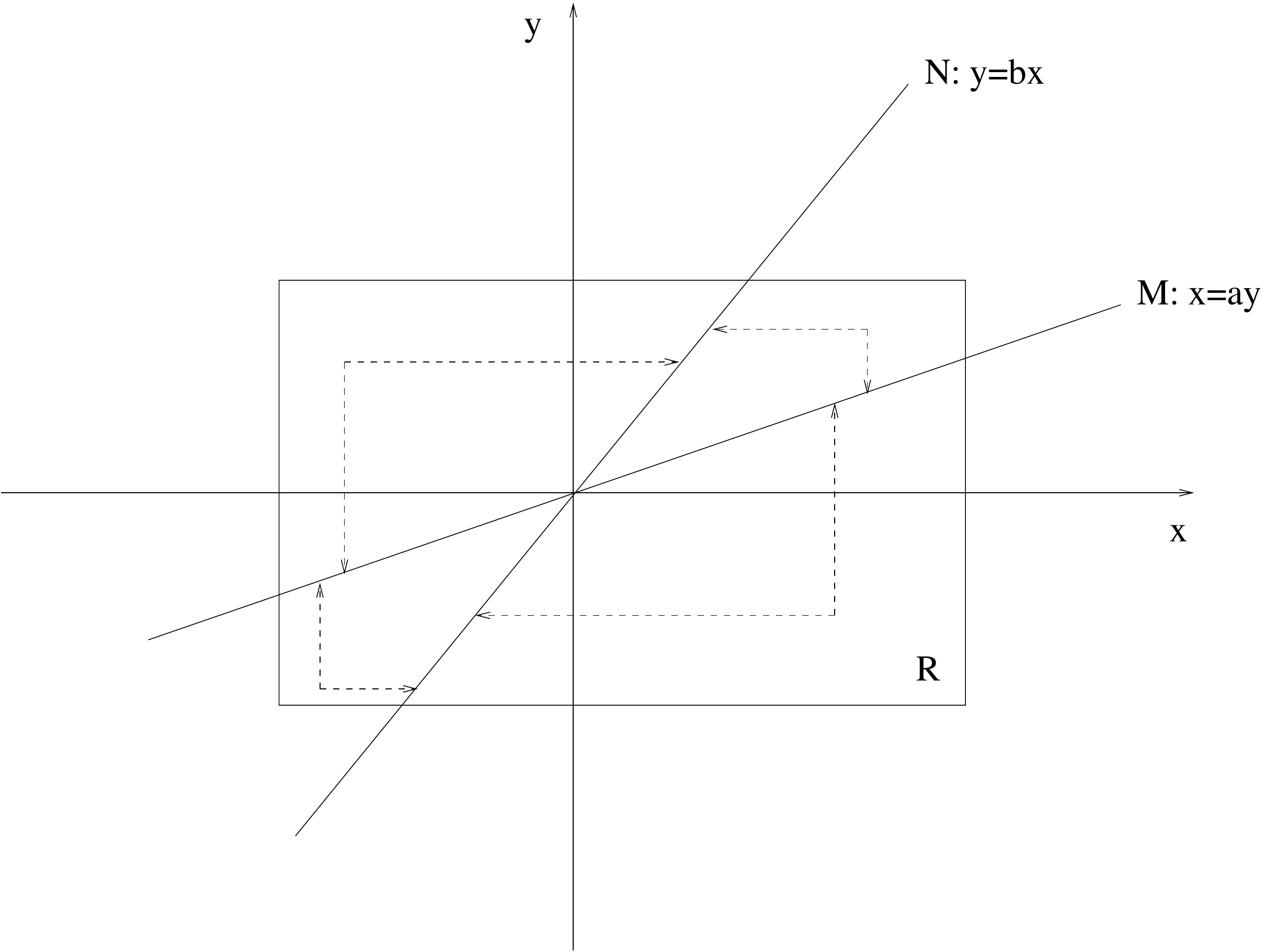}
	\caption{Stable configuration.}
\end{figure} 
%
\begin{figure}\label{Figure_2}
	\centering
	\includegraphics[width=11cm,height=8cm]{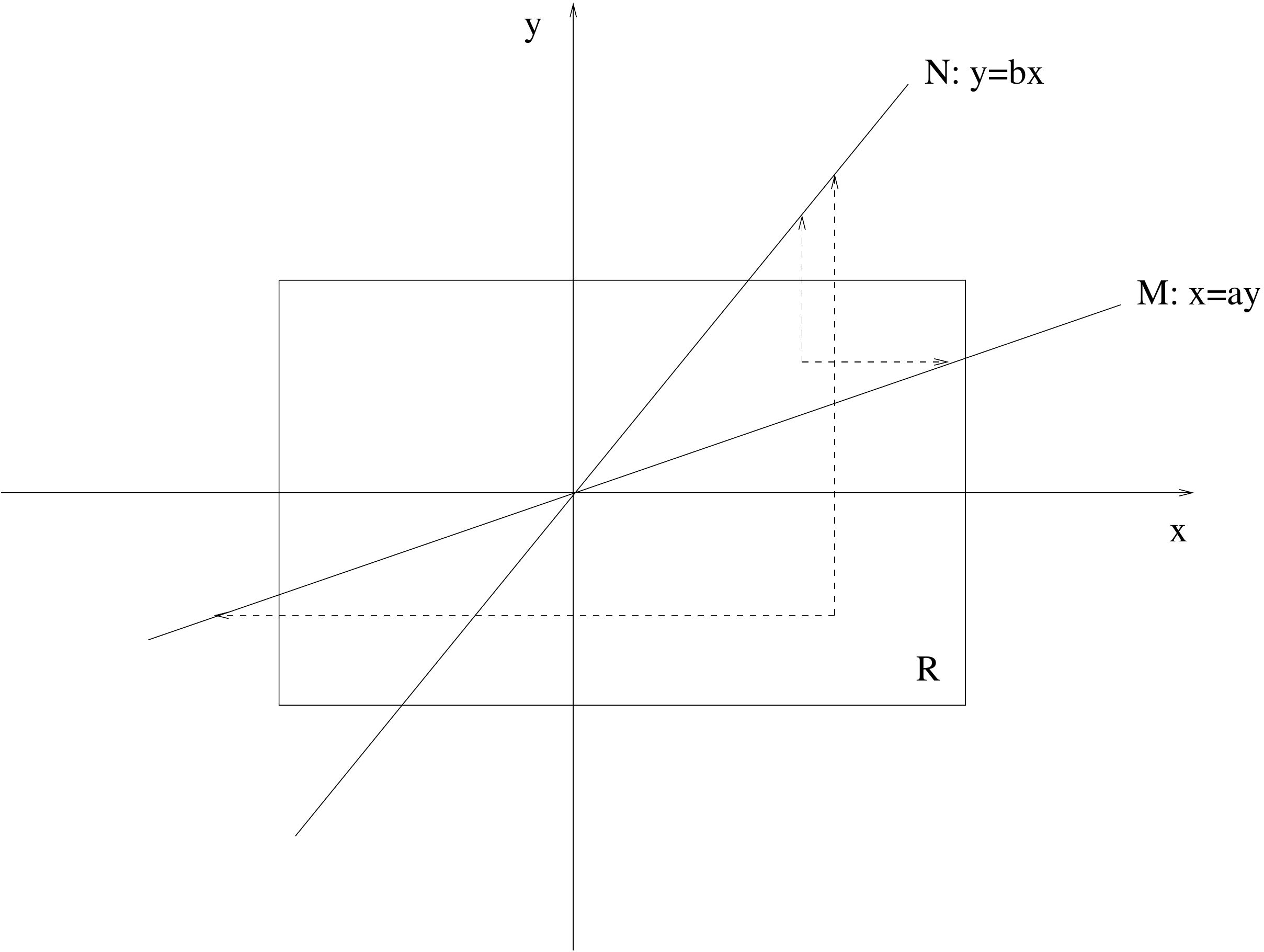}
	\caption{Not a stable configuration.}
\end{figure} 
Clearly, in both cases the transversality condition in (ii) above is met. 
Letting $g$ and $h$ be scalar functions of a single argument
and defined near zero, we now have:
\begin{itemize}
	\item[(a)] With $u(x,b x)=g(x)$  and $v(a y,y)=h(y)$, say, the 
	natural iteration scheme is to set:
	\[u^{(0)}(x,y):=g(x),\qquad v^{(0)}(x,y):=h(y),\]
	and then define
	\[u^{(n+1)}(x,y):=g({\textstyle\frac{y}{b}})+\int_\frac{y}{b}^x v^{(n)}(s,y)\, ds,\]
	and
	\[v^{(n+1)}(x,y):=h({\textstyle\frac{x}{a}})+\int_\frac{x}{a}^y u^{(n)}(x,s)\, ds,\]
	for $n\geq 0$.
	Consider any rectangular neighborhood $R$ of the origin with the 
	property that its upper-right corner and its lower-left corner both lie 
	between the lines $M$ and $N$ in the first and third quadrants, 
	respectively; see Figure 1. It is assumed that $R$ is sufficiently
	small so that $u^{(0)}$ and $v^{(0)}$ are both defined on $R$. 
	It is immediate to verify that any such neighborhood provides a {\em stable} 
	configuration in the following sense: whenever we start at a point $(x_0,y_0)\in R$ 
	and move toward $N$ or $M$ along the integral curves of  
	$\br_1=\del_x$ or $\br_2=\del_y$, respectively, we remain within 
	$R$ until we meet $N$ and $M$, respectively. 
	
	In particular, this guarantees that all iterates $u^{(n)}$ and $v^{(n)}$, 
	as defined above, are well-defined on $R$. We express this by saying that the Stable 
	Configuration Condition (SCC) is satisfied in this case, and this is 
	the situation covered by our theorem below. 
	\item[(b)] The situation changes if we instead prescribe $u$ along $M$ 
	and $v$ along $N$, say, $u(ay,y)=g(y)$ and $v(x,bx)=h(x)$. The 
	natural iteration scheme is now to set:
	\[u^{(0)}(x,y):=g(y),\qquad v^{(0)}(x,y):=h(x),\]
	and then define
	\[u^{(n+1)}(x,y):=g(ay)-\int_x^{ay} v^{(n)}(s,y)\, ds,\]
	and
	\[v^{(n+1)}(x,y):=h(bx)-\int_y^{bx} u^{(n)}(x,s)\, ds,\]
	for $n\geq 0$. The trouble with this is that 
	there is no bounded neighborhood of the origin
	which provides a stable configuration in the above sense; see Figure 2. 
	More precisely, one can verify that given any bounded neighborhood $U$ of the origin, 
	there is always a point $(x_0,y_0)\in U$, such that, either in moving horizontally till 
	intersecting $M$, or in moving vertically till intersecting $N$, one 
	leaves $U$. The upshot is that there is no {\em fixed},
	bounded neighborhood on which all iterates are defined; the SCC is not met
	in this case. 
\end{itemize}
Of course, given locally defined functions $g$ and $h$, we could extend 
them to all of $\RR$, and then run the iteration defined above on all of $\RR^2$. 
However, even assuming convergence to a solution of the original system, 
it appears that the limiting solution would depend on the choice of extensions. 

Indeed, a simple example illustrates that the uniqueness part of our Theorem 
might fail when the SCC is not met. To describe this accurately, we introduce 
the following terminology:
\begin{itemize}
	\item we say that {\em strong uniqueness} holds if there is a fixed 
	neighborhood $U$ of the point $\bar x$ with the property that any two solutions of 
	the system, both of which are defined on $U$, must agree on $U$;
	\item we say that {\em weak uniqueness} holds if, given any two solutions of 
	the system, each of which is defined on an open set about the point $\bar x$, there 
	is a neighborhood $V$ (possibly depending on the given solutions) on
	which the two solutions agree.
\end{itemize}

The uniqueness claim in the Theorem below is that strong uniqueness holds
under the SCC. In contrast, the following simple example shows that 
strong uniqueness might fail when the SCC is not met.
Consider now the trivial system $u_x = 0$, $v_y = 0$, with
constantly vanishing initial data prescribed along the lines $M$ and $N$ depicted in
Figure 2. Consider the system as defined on any open ball about
the origin. For any neighborhood $U$ of the origin on which we wish to define a solution,
since the SCC is not met, there is either 
\begin{itemize}
    \item an interval $(y_0, y_1)$ such that every horizontal line $y =
    \overline y$ with $\overline y \in (y_0,y_1)$ intersects $U$ but only
    meets $M$ outside of $U$, or
    \item an interval $(x_0, x_1)$ such that every vertical line $x = \overline x$ 
    with $\overline x \in (x_0,x_1)$ intersects
    $U$ but only meets $N$ outside of $U$.
\end{itemize}
Suppose the former. Then let $u=F(y)$, where $F$ is a bump function whose
bump is located in $(y_0,y_1)$, and let $v=0$. These functions are a
solution to the system on $U$. Since $u=0,v=0$ is clearly also a solution, we
do not have uniqueness on $U$. (If we instead have the latter situation, we instead
let $v=G(x)$ be a bump function.)

\begin{remark}
	We note that weak uniqueness does hold for the example just considered, 
	even though the SCC is not met. It would be of interest to know if this 
	is the case for more general systems. 
\end{remark}

\begin{remark}
	Returning to the system \eq{u}-\eq{v} above, we note that 
	applying $\del_y$ to the first equation 
	and using the second equation, yield the second order hyperbolic equation 
	$u_{xy}=u$. Conversely, a solution of the latter equation yields, upon
	setting $v:=u_x$, a solution to the system \eq{u}-\eq{v}. Thus, at least 
	at the level of $C^2$-solutions, these are equivalent problems.
	
	Equations of the form $u_{xy}=F(x,y,u,u_x,u_y)$ have been studied 
	extensively, starting with classical treatments by Riemann, Darboux, and Goursat. 
	In particular, various types of boundary value problems
	have been considered; see \cites{khar,lieb} and references therein. 
	However, we are not aware of results that cover the particular situation 
	above, i.e.\ with general data for $u$ and $u_x$ prescribed along two different 
	non-characteristic curves through the origin, and a solution is 
	sought on a full neighborhood of the origin.
\end{remark}

\section{Statement and Proof} \label{sec:theorem}
We start by stating and proving a Lemma about certain restricted systems;
and then state and prove our main theorem. 

We are given open sets $\Omega \subseteq \Rn$ and
$\Upsilon \subseteq \Rm$ together with a $C^1$ frame $\rset_{i=1}^n$ on $\Omega$, 
and we fix a point $\bar x\in \Omega$.
For all positive integers $p$, on all open subsets of $\RR^p$, we will use the max norm, denoted by $\|\cdot \|_\infty$.
The corresponding open ball of radius $\epsilon>0$ about a point $y\in\RR^p$ is denoted $B_\epsilon(y)$.

Throughout this section we use the following conventions:  the integers  $i$, $j$, and $k$  satisfy 
$1 \le i,j,k \le n$ and index  the vector fields in the frame (derivations). 
The integers $\alpha$ and $\beta$  satisfy $1 \le
\alpha,\beta \le m$ and index the unknown functions. 
For each index $\alpha$, $I_\alpha$ will denote the set the indices 
of vector fields with respect to which the unknown function $u_\alpha$ 
is differentiated in the original system \eq{gen_eqns}.
Both the elements and the cardinality of $I_\alpha$ may vary with 
$\alpha$. However, in order to avoid an extra index, whenever $\alpha$ is fixed, we 
simply write $I_\alpha = \{i_1,i_2,\ldots, i_\pa\}$, where it is assumed that
$1\leq i_1 <i_2< \cdots < i_\pa\leq n$.

Finally, for each $\alpha$ we fix an $(n-\pa)$-dimensional $C^1$ submanifold  
$\Xi_\alpha$ of $\Omega$ that is transverse to the span of vector fields 
$\{\br_i\,|\,i\in I_\alpha\}$.  We assume that  the point $\bar x$ belongs to
$\cap_{\alpha=1}^m \Xi_\alpha$ and that each $ \Xi_\alpha$ is small enough 
to be covered by a single coordinate chart centered at $\bar x$. In other words, 
there exist $C^1$-diffeomorphisms $\xi_\alpha\colon \U_\alpha \to \Xi_\alpha$,
where $\U_\alpha \subset\RR^{n-\pa}$ is an open neighborhood of the  origin 
and $\xi_\alpha(0)=\bar x$. The submanifold $\Xi_\alpha$ will be where the 
the unknown function $u_\alpha$ is prescribed. We shall refer to the $\Xi_\alpha$ 
as {\em data manifolds}.

To formulate and prove our key lemma,
  we introduce a collection of new local coordinate charts 
near $\bar x$, one for each $\alpha$. 
For this let $W_i^t$ denote the flow of $\br_i$:
\[\frac{d}{dt}W_i^t(x)=\br_i\big|_{W_i^t(x)},\qquad \qquad W_i^0(x)=x\qquad(i=1,\dots,n).\]
As the vector fields $\br_i$ are $C^1$-smooth, their flows are defined 
and  $C^1$-smooth with respect to $(t,x)$ on an 
$\RR\times\RR^n$-neighborhood of $(0,\bar x)$ (see Theorem 2.6 on 
p.\ 81 in \cite{lang}).

For each $\alpha$, with $I_\alpha = \{i_1,i_2,\ldots, i_\pa\}$, there exists an open neighborhood  of the origin  $\Theta_\alpha\subset \RR^n$
 small enough so that  for each $t=(t_1,\dots, t_n)\in \Theta_\alpha$ the map  
$\psi_\alpha:\Theta_{\alpha} \to \RR^n$ given by
\beq\label{diffeo}
	\psi_\alpha(t):=
	W_{i_{\pa}}^{t_{\pa}} \cdots W_{i_2}^{t_2} W_{i_1}^{t_1}\, \xi_\alpha(t_{\pa+1},\dots,t_n).
\eeq
is well-defined and $C^1$-smooth. Since, by assumption, the manifold 
$\Xi_\alpha$ is everywhere transverse to the span of vector fields 
$\{\br_i\,|\,i\in I_\alpha\}$,
it follows from the $C^1$-version of the inverse mapping theorem 
(Theorem 5.2 on p.\ 13 in \cite{lang})
that, after possibly shrinking $\Theta_\alpha$ near the origin, we may assume that $\psi_\alpha$ is a $C^1$-diffeomorphism between $\Theta_{\alpha}$
and its image $
	\Omega_\alpha:=\psi_\alpha(\Theta_{\alpha})\subset\RR^n.
$
 Thus we now have $m$ coordinate charts    $(\Omega_{\alpha}, \psi_\alpha^{-1})$ near $\bar x$.
 From now on we shrink $\Xi_\alpha$ to $\Xi_\alpha\cap\Omega_\alpha$.

Lemma~\ref{lemma} below can be regarded as a ``restricted'' form of our main theorem, in
which integrability conditions are not assumed and, consequently, we cannot
conclude that the equations in the original system \eq{gen_eqns} have solutions 
defined on a full $\RR^n$-neighborhood of $\bar x$. Instead, we
will only conclude that each equation is satisfied on a certain submanifold
of $\Omega$ (which varies from equation to equation).
To define these submanifolds, for each $\alpha$,
we introduce a sequence of manifolds 
$\Xi_\alpha^0,\dots,\Xi_\alpha^\pa$  by:
\beq\label{Xij}\Xi_\alpha^j:=\{\psi_\alpha(t)\,|\, t=(t_1,\dots,t_j,0,\dots,0,t_{\pa+1},\dots,t_n)\in\Theta_\alpha \}.\eeq
We observe that $\Xi_\alpha^0\equiv\Xi_\alpha$,  
 while  $\Xi_\alpha^j$ is the set of points obtained by starting from a 
point  in $\Xi_\alpha$ and then applying the flows $W_{i_1}^{t_1},\dots,W_{i_j}^{t_j}$,
in that order.  In particular, 
$\Xi_\alpha^\pa \equiv \Omega_\alpha$.

The existence of the restricted solutions in Lemma~\ref{lemma} will be obtained via Picard iteration.
As outlined in Section \ref{scc_cond} on the Stable Configuration Condition (SCC),
we shall need to impose constraints on 
the coordinate charts $(\Omega_{\alpha}, \psi_\alpha^{-1})$ in order to guarantee that 
the various iterates are well-defined. To avoid a further level of technical details we 
choose to require that the SCC is satisfied in arbitrarily small neighborhoods. We
first introduce the following terminology. 
\begin{definition}\label{def-access}  
    	Let $p\in\{1,\dots,n\}$. An open neighborhood $\Theta \subset\RR^n$ of the origin is 
    	called \emph{$p$-accessible} if for each $t=( t_1, t_2, \ldots,t_n)\in \Theta$,
    	the piecewise linear path from $ (0,\dots,0, t_{p+1},\dots,t_n)$ to $t$, with 
    	vertices 
	\begin{align*}
	&(0,\dots,0, t_{p+1},\dots,t_n),\\
	&(t_1,0,\dots,0,t_{p+1}, \dots,t_n),\\
	&(t_1, t_2,0\dots,0,t_{p+1},\dots,t_n),\\
	&\qquad\vdots\\
	&(t_1,t_2, \dots,t_n),
	\end{align*}
    	belongs to $\Theta$.
\end{definition}
\begin{definition}[Stable Configuration Condition]\label{def-stable}  
In the context described above, let $\R=\{\br_1,\dots,\br_n\}$ be a
local frame near $\bar x$, and let a set of $m$ submanifolds 
$\Xi:=\{\Xi_1,\dots, \Xi_m\}$, of co-dimensions $p_1,\dots,p_m$ in $\RR^n$, 
respectively, and containing $\bar x$, be given. 
We say that $(\R,\Xi)$ is a {\em stable configuration} near 
$\bar x\in \Omega$ if for every $\epsilon >0$, there exists an 
open neighborhood  $\Omega_\epsilon\subset B_\epsilon(\bar x)$ 
of $\bar x$ and   {$\pa$-accessible} neighborhoods 
$\Theta_{\alpha,\epsilon}\subset \RR^{n}$, $\alpha=1,\dots,m$, 
such that for all $\alpha$  the map  
$\psi_\alpha\colon \Theta_{\alpha,\epsilon} \to\Omega_\epsilon$, 
defined by \eqref{diffeo} is a $C^1$-diffeomorphism. 
The neighborhoods $\Omega_\epsilon$ are called \emph{stable neighborhoods}.
\end{definition}

\noindent The significant part of this definition is that the same neighborhood $\Omega_\epsilon$
works for all $\alpha=1,\dots,m$. Also, $\Omega_\epsilon$ is stable (or invariant) in the 
sense that, given any point $x\in\Omega_\epsilon$, we can flow back to each $\Xi_\alpha$
along integral curves of the vector fields $\{\br_i\,|\,i\in\Ia\}$ (in decreasing order) without leaving 
$\Omega_\epsilon$. The following lemma provides a key technical step and will be applied twice 
in the proof of the main theorem.

\begin{lemma}\label{lemma}
    Let $\Omega \subseteq \Rn$ and
$\Upsilon \subseteq \Rm$ be open subsets. For $\alpha=1,\dots, m$, let   $I_\alpha = \{i_1,i_2,\ldots, i_\pa\}$ be an ordered set of indices: 
$1\leq i_1 <i_2< \cdots < i_\pa\leq n $.  Assume the following:
    \begin{enumerate}

        \item $\R=\rset_{i=1}^n$ is a $C^1$ frame on $\Omega$;

        \item for each $\alpha$, $\Xia$ is an embedded $C^1$ manifold
        in $\Omega$, and $\bar x \in \cap_\alpha
        \Xia$;

        \item the manifold $\Xia$ is of codimension $p_\alpha$ and is
        everywhere transverse to the span of $\rset_{i\in I_\alpha}$;

\item the frame $\R$ and the set of manifolds 
$\Xi:=\{\Xi_1,\dots, \Xi_m\}$ is a \emph{stable configuration} near $\bar x\in \Omega$
according to Definition \ref{def-stable};

        \item for each  $\alpha$, $\ga:\Xia \to \RR$ is $C^1$-smooth and bounded, and
        \[\og := (g_1(\ox), g_2(\ox), \ldots, g_m(\ox)) \in \Upsilon;\]
  \item for each $\alpha$ and each $i\in I_\alpha$, $f_i^\alpha:\Omega
        \times \Upsilon \to \RR$ is uniformly bounded and continuous on $\Omega
        \times \Upsilon$ and also uniformly Lipschitz in the second argument.

    \end{enumerate}
    Then there is a neighborhood $\wtO \ni \ox$ on which there is a unique
    solution $u:\wtO\to\Upsilon$ to the system
    \begin{equation} \label{eq:restricted}
        \br_{i_j}(\ua)\atx = f^\alpha_{i_j}(x, u(x)) \quad
        \text{for $1 \le \alpha \le m$, $i_j \in \Ia$,
       and $x \in \Xia^j \cap \wtO$}
    \end{equation}
    satisfying the data
    \begin{equation} \label{eq:data}
        \ua(x) = \ga(x) \; \textrm{ for } x \in \Xia \cap \wtO .
    \end{equation}
\end{lemma}
%
\begin{remark}\label{doubl_indx}
	Note that, in \eq{eq:restricted}, we need to employ double subscripts $\ijj$ 
	for the elements of the ordered set $I_\alpha$ since the position (the $j$th, say) of 
	an index in $I_\alpha$ relates to which submanifold $\Xia^j$ is considered.
\end{remark}
\begin{proof} 
Let $L$ be a common Lipschitz constant for the functions $\fia$:
    \beq\label{Lip}
		|f_i^\alpha(x,u)-f_i^\alpha(x,v)|\leq L\|u-v\|_\infty
		\qquad\text{whenever $x\in \Omega$, $u,v\in\Upsilon$.}
	\eeq
We let $M$ be a common bound for the functions $\fia$:
    \beq \label{M} |\fia(x,u)| \le M
    \text{ for $(x,u) \in \Omega\times \Upsilon$.}
    \eeq
    Choose $r > 0$ such that $\bar B_r(\og)$ (the closed ball under the sup norm)
    is contained in $\Upsilon$ and put
    \beq\label{r}
    	\wtU := \bar B_r(\og).
    \eeq
Shrink $\Omega$ to $\Omega' \ni \ox$ such that
\beq\label{g}
   	|\ga(x) - \ga(\ox)| \le {\textstyle\frac{1}{2}}r
   	\quad\text{for each $\alpha$ and
    	for $x \in \Xia \cap \Omega'$.}
\eeq
    Finally,  using assumption (4), we choose a bounded and stable neighborhood 
    $\wtO \subset \Omega'$, with accessible neighborhoods $\Theta_\alpha$ 
    and diffeomorphisms  $\psi_\alpha\colon \Theta_{\alpha} \to\wtO$, according to 
    Definitions~\ref{def-access},~\ref{def-stable}. If necessary, we may shrink $\wtO$, and hence
    each $\Theta_\alpha$, so as to have
\beq\label{Q}
	nL \|t\|_\infty \le {\textstyle\frac{1}{2}} \quad \text{and}\quad 
	nM \|t\|_\infty \le {\textstyle\frac{1}{2}}r
    	\quad  \text{for all $t\in \Theta_\alpha$, $\alpha=1,\dots,m$.}
\eeq  
Next, let $\mC$ denote the set of continuous functions $\wtO \to \wtU$. Define a functional
    $\Phi:\mC \to \mC$ by defining its components according to
    \begin{align} 
     &\Phi[u]_\alpha(x) := g_\alpha(\xia)
        + \int_0^{t_1} f_{i_1}^\alpha(
            \psia( s, 0,  \ldots, 0, t_{\pa+}),
            u(\psia( s, 0,  \ldots, 0, t_{\pa+}))
        ) \, ds\nn \\
        &+ \int_0^{t_2} f_{i_2}^\alpha(
            \psia( t_1,s, 0,  \ldots, 0, t_{\pa+}),
            u(\psia( t_1,s, 0,  \ldots, 0, t_{\pa+}))
            ) \, ds\nn \\
        &\qquad \vdots\label{eq:functional} \\
        &+ \int_0^{t_\pa} f_{i_\pa}^\alpha(
            \psia( t_1, \ldots, t_{\pa-1}, s, {t_{\pa+}}),
            u(\psia( t_1, \ldots, t_{\pa-1}, s, {t_{\pa+}}))
            ) \, ds .\nn
    \end{align}
    In the above equation, the values $(t_1, t_2, \ldots, t_n)$ are
    chosen so that $x = \psia( t_1,  \ldots, t_n)$. Since $\psia$, given by \eq{diffeo}
    is a diffeomorphism, this choice is unique.  We have also used the abbreviations:  
    $t_{\pa+}:=(t_{\pa+1},\dots, t_n)$, $\xia:=\xia(t_{\pa+})$. 
   Note that the function $\Phi[u]$ is well-defined since $u\in \mC$ 
   and since the neighborhoods $\Theta_\alpha$ are accessible; this is the
   technical reason for imposing the Stable Configuration Condition.

    To verify that $\Phi$ in fact maps $\mC$ into itself, we assume $u \in \mC$ and 
    $x \in \wtO$ and show that the right hand side of
    $\eqref{eq:functional}$ belongs to $\wtU = \bar B_r(\og)$. To this end, observe that
    \newcommand{\myphantom}{
        \left| \sum_{j=1}^\pa \int_0^{t_j} 
        f_{\ia^j}^\alpha(\psia, u(\psia)) \, dx
        \right|
    }
    \begin{align}
        & \| \Phi[u](x) - \og \|_\infty
       = \max_\alpha \left| \Phi[u]_\alpha(x) - \ga(\ox) \right| \nn\\
        \label{eq:inc-1}
        &\qquad\qquad\le \max_\alpha \Big\{ \vphantom{\myphantom} | \ga(\xi) - \ga(\ox) | +
            \Big| \sum_{j=1}^\pa \int_0^{t_j} 
            f_{\ia^j}^\alpha(\psia, u(\psia)) \, ds
            \Big| \Big\} \\
        \label{eq:inc-2}
        &\qquad\qquad\le \max_\alpha \Big( {\textstyle\frac{1}{2}}r + n M\cdot \big(\max_{1\leq j\leq \pa} |t_j|\big) \Big)\le r ,
    \end{align}
    where we have omitted the arguments of $\psi_\alpha$. Here, to obtain line \eqref{eq:inc-1} we used the definition
    \eqref{eq:functional} of the functional $\Phi$ and the triangle
    inequality. To obtain \eqref{eq:inc-2} we used statement \eq{g} for the first term and
 the triangle inequality     together with  \eq{M} and the fact
    that $\pa \le n$ for the second term. Finally, we have used \eq{Q}.

    Equip $\mC$ with the metric $d(u,v) = \sup_{x\in\wtO}
    \|u(x)-v(x)\|_\infty$. With this metric, $\mC$ is a complete metric space.
    We now show that $\Phi$ is a contraction mapping. Let $u,v \in \mC$, and estimate
    \begin{align}
        \nonumber
        d(\Phi[u], \Phi[v])
        &= \sup_{x\in\wtO} \|\Phi[u](x) - \Phi[v](x)\|_\infty = \sup_{x\in\wtO} \max_\alpha
            \left|\Phi[u]_\alpha(x) - \Phi[v]_\alpha(x) \right| \nn\\
        \label{eq:cont-1}
        &\le \sup_{x\in\wtO} \max_\alpha
            \sum_{j=1}^\pa \int_0^{t_j^\alpha} \left|
             f_{i_j}^\alpha(\psia, u(\psia)) -
             f_{i_j}^\alpha(\psia, v(\psia))\right| \, ds \\
        \label{eq:cont-2}
        &\le n  \max_{{1\leq \alpha\leq m}\atop{ 1\leq  j\leq \pa}} | t_j^\alpha |\,L\cdot\sup_{y\in\wtO}  \|u(y) - v(y)\|_\infty \\
        \label{eq:cont-3}
        &\le {\textstyle\frac{1}{2}}d(u, v),
    \end{align}
    where, again, we have omitted the arguments of $\psia$. To obtain line
    \eqref{eq:cont-1}, we used the definition \eqref{eq:functional} of $\Phi$
    and the triangle inequality. To obtain line \eqref{eq:cont-2}, we used
    the fact that each $\fia$ has Lipschitz constant $L$ and that $p_\alpha
    \le n$. To obtain line \eqref{eq:cont-3}, we used  the first inequality in \eq{Q} and the definition 
    of $d(u,v)$. Thus $\Phi$ is a (uniformly) strict contraction. It follows that $\Phi$ 
    has a unique fixed point in $\mC$, which we denote $\wu$.  Here we used notation $t_j^\alpha$, to emphasize  that the pre-image of $x$ under $\psi_\alpha$ depends on $\alpha$. Thus
    \begin{align} 
        &\wu_\alpha(x)  = g_\alpha(\xia)
        + \int_0^{t_1} f_{i_1}^\alpha(
            \psia( s, 0,  \ldots, 0, t_{\pa+}),
            \wu(\psia( s, 0,  \ldots, 0, t_{\pa+}))
        ) \, ds\nn \\
        &+ \int_0^{t_2} f_{i_2}^\alpha(
            \psia( t_1,s, 0,  \ldots, 0, t_{\pa+}),
            \wu(\psia( t_1,s, 0,  \ldots, 0, t_{\pa+}))
            ) \, ds\nn \\
        &\qquad \vdots\label{eq:fixed-point} \\
        &+ \int_0^{t_\pa} f_{i_\pa}^\alpha(
            \psia( t_1, t_2, \ldots, t_{\pa-1}, s, {t_{\pa+}}),
            \wu(\psia( t_1, t_2, \ldots, t_{\pa-1}, s, {t_{\pa+}}))
            ) \, ds .\nn   
    \end{align}
    Since $t_1=\cdots=t_\pa=0$ and $\xia\equiv\xia(t_{\pa+})=x$, whenever $x\in \Xia$, the function $\wu$ satisfies the data \eq{eq:data}.
    Note that on the manifold $\Xia^j$, we have $t_k = 0$ for $j<k\leq \pa$, and
    also $\br_{i_j} = \partial_{t_j}$. The latter follows since, for any
    smooth function $h: \wtO \to \RR$, whenever $x \in \Xia^j$, we have
    \beq\label{tj}
        \partial_{t_j} h(x)
        = \partial_{t_j} h(W_j^{t_j} \cdots W_1^{t_1} \xia) 
        = \nabla_xh \cdot \br_{i_j}\big|_{W_j^{t_j} \cdots W_1^{t_1} \xia} 
        = \br_{i_j} h(x) .
    \eeq
    Thus, with equation \eqref{eq:fixed-point} restricted to
    $\Xia^j$, we apply the fundamental theorem of calculus and obtain that, for $x\in \Xia^j$,
    \begin{align*}
      &  \br_{i_j} \wu_\alpha(x)
        = \partial_{t_j} \wu_\alpha(x) \\
        &\quad= f_{i_j}^\alpha(\psia(t_1, \ldots, t_j, 0, \ldots, 0, t_{\pa+1}),
        \wu(\psia(t_1, \ldots, t_j, 0, \ldots, 0, t_{\pa+1})))\\
       & \quad= f_{i_j}^\alpha(x, \wu(x)),
    \end{align*}
    showing that $\wu$ is indeed a solution of \eq{eq:restricted}-\eq{eq:data}.

    It remains to show that $\wu$ is the unique solution 
    of \eq{eq:restricted}-\eq{eq:data} on $\wtO$. 
    Assuming $v:\wtO\to\Upsilon$ is also a solution to \eq{eq:restricted}-\eq{eq:data} in $\wtO$, we have
    \begin{align}
        d(\wu,v)&=\sup_{x \in \wtO} \| \wu(x)- v(x) \|_\infty 
        = \sup_{x \in \wtO} \max_\alpha |\wu_\alpha(x) - v_\alpha(x)|\nn \\
        \label{eq:unique2}
        &=\sup_{x \in \wtO} \max_\alpha
           \Big| \sum_{j=1}^{p_\alpha} \int_0^{t_{i_j}}
           f_{i_j}^\alpha(\psia, \wu(\psia)) - f_{i_j}^\alpha(\psi_\alpha, v(\psi_\alpha)) dx \Big| \\
        \label{eq:unique3}
        &\le \sup_{y \in \wtO} \max_\alpha
          \sum_{j=1}^{\pa} |t_{i_j}| L \| \wu(y) - v(y) \|_\infty \\
        \label{eq:unique4}
        &\le {\textstyle\frac{1}{2}}\sup_{y \in \wtO} \|\wu(y) - v(y)\|_\infty ={\textstyle\frac{1}{2}}d(\wu,v),
    \end{align}
    where, as above, $x = \psia( t_1,  \ldots, t_n)$ and we have omitted the arguments of $\psia$. 
    Line \eqref{eq:unique2}
    follows from the fact that both $\wu$ and $v$, being solutions to \eq{eq:restricted}-\eq{eq:data}, satisfy
    \eqref{eq:fixed-point}. Line \eqref{eq:unique3} follows from the triangle
    inequality, the Lipschitz property \eq{Lip}, and the fact that 
    the neighborhoods $\Theta_\alpha$ are accessible (so that the $\psi_\alpha$ 
    take values in $\wtO$). The inequality in  \eqref{eq:unique4}
    follows from the first inequality in \eq{Q}. Now, $\tilde \Omega$ is, by choice, a bounded set
    in which each point can be reached by starting on any data manifold $\Xi_\alpha$
    and then moving along a finite number of integral curves of the vector fields $\br_i$. 
    Since we assume that the initial data $g_\alpha$ are bounded, and also that the functions 
    $f_i^\alpha:\Omega\times \Upsilon \to \RR$ are uniformly bounded, it follows that 
    the solutions $\tilde u$ and $v$ are both bounded on $\tilde \Omega$.
    Therefore, $\sup_{x \in \wtO} \| \wu(x) - v(x) \|_\infty < \infty$ and 
    it follows from the inequalities above that $d(\wu,v)=0$, i.e.\ $\wu\equiv v$ on $\wtO$.
\end{proof}

\begin{theorem}\label{thm}
    Suppose that, in addition to hypotheses (1)--(6) of  Lemma~\ref{lemma}, we also have that:
    \begin{enumerate}

        \item[(7)] the functions $\fia$ belong to $C^1(\Omega\times \Upsilon)$;

        \item[(8)] for each $\alpha$, for each $i,j \in \Ia$ with $i \ne j$,
        and for each $\beta$: if $i\notin I_\beta$, then $\partial_{u_\beta}
        f_j^\alpha = 0$;

        \item[(9)] for each $\alpha$, the vector
        fields $\rset_{i \in \Ia}$ are in involution, i.e.  $[r_j,r_k]\in \spa\rset_{i \in \Ia}$ 
        whenever $j,k\in \Ia$;

        \item[(10)] for each $\alpha$, for each $i,j \in \Ia$ with $i \ne j$, and for 
        all $(x,u)\in\Omega\times\Upsilon$:
        \begin{align} \label{eq:integrability}
            &(\nabla_x f_j^\alpha)|_{(x,u)} \cdot \br_i\atx
            + \sum_{\beta: i \in I_\beta} \partial_{u_\beta} f_j^\alpha(x,u) f_i^\beta(x,u) \\
            \nonumber &- \nabla_x(f_i^\alpha)|_{(x,u)} \cdot \br_j\atx
            - \sum_{\beta: j \in I_\beta} \partial_{u_\beta} f_i^\alpha(x,u) f_j^\beta(x,u) =
            \sum_{k \in \Ia} c_{ij}^k(x) f_k^\alpha(x, u),
        \end{align}
	where here and below $\nabla_x(f^\alpha_i)$  denotes the gradient 
	with respect to the variables $x$ and $c_{ij}^k$ denote the structure 
	coefficients of the frame:
	\[[\br_i,\br_j]\big|_x = \sum_{k=1}^n c_{ij}^k(x)\br_k\big|_x
	\qquad \text{for $x\in \Omega$ and $1\leq i,j\leq n$}.\]
        
        \item[(11)] for each $\alpha$ and each $i,j,k\in \Ia$, the structure coefficient $c_{ij}^k(x)$
        is uniformely bounded on $\Omega$.
\end{enumerate}
    Then there is a neighborhood $\wtO$ of $\ox$ on which there is a unique solution to
    the system
    \begin{equation} \label{eq:system}
        \br_{i}(\ua)\atx = \fia(x, u(x)) \quad\text{for $1 \le \alpha \le m$, $i \in \Ia$, and $x \in \wtO$,}
    \end{equation}
    satisfying the data
    \begin{equation} \label{eq:data2}
        \ua(x) = \ga(x) \; \textrm{ for } x \in \Xia \cap \wtO .
    \end{equation}
\end{theorem}
Before giving the proof we make a few remarks.
First, as remarked above, the difference between the conclusion of Theorem~\ref{thm} and
the conclusion of Lemma~\ref{lemma} is that in the Theorem, the equations of system
\eqref{eq:system} are satisfied everywhere in $\wtO$, while in the Lemma,
each equation of system \eqref{eq:restricted} is only guaranteed to be
satisfied only for   $x\in \Xia^j \cap \wtO$. The data $\eqref{eq:data}$ and
$\eqref{eq:data2}$ are identical.

The integrability conditions appearing in equations \eqref{eq:integrability}
are the generalization of the condition of mixed partial derivatives being equal,
to the case of non-commutative derivations.
They correspond to the integrability conditions in the
PDE version of the classic Frobenius Theorem
(see the equations marked (**) in Theorem 1 of Chapter 6 in \cite{spi}).
The conditions appear from expanding the equation
\begin{equation} \label{eq:prelude-integrability}
    \br_i(\br_j(u_\alpha)) - \br_j(\br_i(u_\alpha)) = \sum_{k=1}^n c_{ij}^k
    \br_k(u_\alpha),
\end{equation}
which should hold for any function $u_\alpha$, and then, once fully expanded,
making substitutions of the form $\br_i(u_\alpha) = f_i^\alpha(x,u)$, which
should hold for any solution $u = (u_1, \ldots, u_m)$ of
\eqref{eq:system}.

The restricted summations of the form $\{\beta: i \in I_\beta\}$,
$\{\beta: j \in I_\beta\}$, and $\{k \in I_\alpha\}$ in \eqref{eq:integrability} 
ensure that \eqref{eq:integrability} only contains functions $f_i^\alpha$ which actually
are defined by the system \eqref{eq:system}. For instance, examining the
first summation in \eqref{eq:integrability}, if we included an index $\beta$
with $i \notin I_\beta$, then the factor $f_i^\beta(x,u)$ would be a function
that is not given by the system \eqref{eq:system}. Similar remarks apply to
the other two summations.

Hypotheses (8) and (9) are necessary so that in making the restrictions on
summation indices just described, we have not actually omitted any terms that
should appear in the expansion of \eqref{eq:prelude-integrability}.
For instance, for the omitted indices $\beta$ with $i \notin I_\beta$ from
the first summation of \eqref{eq:integrability}, hypothesis (8) guarantees that
$\partial_{u_\beta} f_j^\alpha(x,u) \equiv 0$, and so we have not actually
missed any terms.
Similar remarks apply to hypothesis (9) and the last summation.

\begin{proof}
    Apply  Lemma~\ref{lemma} to obtain a neighborhood $\wtO$ of $\ox$ on which
    there is a unique solution $\wu$ to the system \eqref{eq:restricted}
    satisfying \eqref{eq:data2}. It remains only to show that $\wu$ is a
    solution to the full system \eqref{eq:system}; uniqueness is already
    established since any solution of \eqref{eq:system} is also a solution of
    \eqref{eq:restricted}.

    Fix (for now) an index $\alpha$. As in the formulation of Lemma \ref{lemma},
    we will use double indices $\ijj$ for the elements of the ordered set 
    $I_\alpha=\{i_1,\dots,i_\pa\}$.
    For each $\ijj\in\Ia$ (i.e., $1\leq j\leq \pa$) we define the function
    $A_{\ijj}^\alpha:\wtO\to\RR$ by
    \beq\label{A}
        A_{i_j}^\alpha(x) := \br_{i_j}(\tilde u_\alpha)\atx - f^\alpha_{i_j}(x,\wu(x)).
    \eeq
    Then, for each $i_k\in\Ia$ with $j < k\leq \pa$, we apply $\br_{i_k}$ to \eq{A} to obtain 
    \begin{align}
        \br_{i_k}(A_{i_j}^\alpha)\atx &=
        \br_\ikk(\br_\ijj(\wua))\atx - (\nabla_x \fa_\ijj)|_{(x, \wu(x))}\cdot\br_\ikk\atx \nn\\
        &\quad \quad - \sum_{\beta: \ikk \in I_\beta} \partial_{u_\beta} \fa_{\ijj}(x,\wu(x)) \br_\ikk(\wu_\beta)\atx 
        \label{eq:special1}\\
        \label{eq:special2}
        &= \br_\ijj(\br_\ikk(\tilde u_\alpha))\atx + \sum_{l\in\Ia} c_{\ikk \ijj}^l(x) \br_l(\tilde u_\alpha)\atx \\
        \nonumber
        &\quad \quad - (\nabla_x \fa_\ijj)|_{(x, \wu(x))}\cdot\br_\ikk\atx
        - \sum_{\beta: \ikk \in I_\beta} \partial_{u_\beta} \fa_{\ijj}(x,\wu(x)) \br_\ikk(\wu_\beta)\atx.
    \end{align}
    We note that the summation in line \eqref{eq:special1} is restricted to
    $\{\beta: \ikk\in I_\beta\}$ by hypothesis (8), and that the
    summation in line \eqref{eq:special2} is restricted to $\{l \in I_\alpha\}$ by
    hypothesis (9).
    We now restrict the last equation above to $x \in \Xia^k$, where, according to the 
    conclusion of Lemma \ref{lemma}, $\br_\ikk(\tilde u_\alpha) = f_\ikk^\alpha(x,\wu)$. Thus, 
    for $x \in \Xia^k$ we have
    \begin{align}
        \br_{i_k}(A_{i_j}^\alpha)\atx
        &= \br_\ijj(f_\ikk^\alpha(x,\wu))\atx + \sum_{l\in\Ia} c_{\ikk\ijj}^l(x) \br_l(\tilde u_\alpha)\atx \\
        \nonumber
        &\quad  - (\nabla_x\fa_\ijj)|_{(x,\wu)}\cdot\br_\ikk\atx
        - \sum_{\beta:\ikk\in I_\beta} \partial_\beta \fa_\ijj(x,\wu) \br_\ikk(\wu_\beta)\atx \\
        \label{eq:about-to-subtract}
        &= (\nabla_xf_\ikk^\alpha)|_{(x,\wu)} \cdot \br_\ijj\atx
        + \sum_{\beta: \ikk \in I_\beta} \partial_{u_\beta} f_\ikk^\alpha(x,\wu) \br_\ijj(\wu_\beta)\atx\nn\\
        &\qquad\qquad\qquad\qquad\qquad\qquad\qquad+ \sum_{l \in I_\alpha} c_{\ikk\ijj}^l(x) \br_l(\tilde u_\alpha)\atx \\
        \nonumber
        &\quad - (\nabla_xf_\ijj^\alpha)|_{(x,\wu)} \cdot \br_\ikk\atx
        - \sum_{\beta: \ikk \in I_\beta} \partial_{u_\beta} f_\ijj^\alpha(x,\wu) \br_\ikk(\wu_\beta)\atx .
    \end{align}
    (Here we begin omitting the argument $x$ of $\wu$). We now apply the integrability condition 
    \eqref{eq:integrability} with $u =\wu\equiv \wu(x)$, $i=i_j$, $j=i_k$, and rearrange to obtain
    \begin{multline} \label{eq:integrability-rearranged}
        (\nabla_x \fa_\ikk)|_{(x,\wu)} \cdot \br_\ijj\atx
        + \sum_{\beta: \ijj \in I_\beta} \partial_{u_\beta} f_\ikk^\alpha(x,\wu) f_\ijj^\beta(x,\wu)
        + \sum_{l\in \Ia} c_{\ikk\ijj}^l(x) f_l^\alpha(x,\wu) \\
        - (\nabla_x f_\ijj^\alpha)|_{(x,\wu)}\cdot \br_\ikk\atx
        - \sum_{\beta: j \in I_\beta} \partial_{u_\beta} f_\ijj^\alpha(x, \wu) f_\ikk^\beta(x,\wu) = 0 .
    \end{multline}
    Subtracting \eqref{eq:integrability-rearranged} from
    \eqref{eq:about-to-subtract} we obtain that: whenever $\ijj, \ikk \in \Ia$ with $k > j$,
    and $x \in \Xia^k$, then
    \begin{align} \label{eq:final-system}
        \br_\ikk(A_\ijj^\alpha)\atx
        &= \sum_{\beta: \ijj \in I_\beta} \partial_{u_\beta} f_\ikk^\alpha(x,\wu) A_\ijj^\beta(x) 
        - \sum_{\beta: \ikk \in I_\beta} \partial_{u_\beta} f_\ijj^\alpha(x,\wu) A_\ikk^\beta(x)\nn\\
        &\quad+ \sum_{l\in \Ia} c_{\ikk\ijj}^l(x) A_l^\alpha(x).
    \end{align}
    We now consider \eqref{eq:final-system} as a system of differential equations for the 
    unknowns $A_\ijj^\alpha$, with vanishing data prescribed along appropriate submanifolds. 
    The trivial functions $A_\ijj^\alpha\equiv 0$ clearly provide a solution.
    The goal is to apply the uniqueness part of Lemma~\ref{lemma} to the new system 
    \eqref{eq:final-system}, and conclude that the trivial solution is the only one. That is, 
    $A_\ijj^\alpha$  defined by \eqref{A} must vanish identically near $\bar x$.
  We proceed with verifying  that 
    the assumptions of Lemma~\ref{lemma} are satisfied for the system \eq{eq:final-system}.
    
    To do so we introduce the following notations. Let 
    $J=\{(j,\alpha)\,|\,1 \le \alpha \le m; \; 1\leq j \leq \pa\}$, and for each double
    index $(j, \alpha)\in J$ define
    $\Lambda_{j, \alpha} := \Xia^j$, where the right-hand side is given by \eq{Xij}.  
    Thus, the $C^1$ submanifolds $\Lambda_{j, \alpha}$ can be  parametrized  
    by the functions
    \begin{align*} 
	\lambda_{j,\alpha}(s_1,\dots, s_{n-\pa+j})&:= W_{\ijj}^{{s_j}} 
	\cdots  W_{i_{1}}^{s_{1}}\, \xi_{\alpha}(s_{j+1},\dots,s_{n-\pa+j})\\
	&\;= \psia(s_1,\dots,s_j,0,\dots,0,s_{j+1},\dots,s_{n-\pa+j}),
     \end{align*}
     which are defined for all $(s_1,\dots, s_{n-\pa+j})$ with the property that the
     point $(s_1,\dots,s_j,0,\dots,0,s_{j+1},\dots,s_{n-\pa+j})$ belongs to $\Theta_\alpha$.
     The domain of $\lambda_{j,\alpha}$ is, therefore, the intersection of $\Theta_\alpha$ with the
     coordinate subset of $\RR^n$ where $t_{j+1}=\cdots=t_{\pa}=0$.

     We let the data $h_{j,\alpha}$ for the unknown $A^\alpha_{i_j}$ in the 
     system \eq{eq:final-system} vanish identically 
     on the manifolds  $\Lambda_{j,\alpha}$: $h_{j,\alpha}(x) \equiv 0$ for 
     $x\in \Lambda_{j,\alpha}$.  Finally, setting $J_{j,\alpha}:= \{ k \,|\, j<k\leq\pa \}$, 
     we have that \eq{eq:final-system} yields an equation for $\br_\ikk(A_\ijj^\alpha)$ 
     whenever $(j,\alpha)\in J$ and $k\in J_{j,\alpha}$.
    
    Now consider the hypotheses of Lemma~\ref{lemma} in the context of  system \eq{eq:final-system}:
    \begin{itemize}

        \item[(1)'] $\R=\rset_{i=1}^n$ is a $C^1$ frame on $\wtO$;

        \item[(2)']  $\Lambda_{j,\alpha}$, $ (j,\alpha)\in J$, are
        embedded $C^1$ submanifolds of $\wtO$, and they share a common point $\ox \in
        \cap_{(j,\alpha)\in J} \Lambda_{j,\alpha}$;

        \item[(3)']  each manifold $\Lambda_{j,\alpha}$ is of codimension $\pa-j$ and is  everywhere transverse to
        the span of $\{\br_\ikk\}_{k \in J_{j,\alpha}}$;
        
        \item[(4)']  \label{h-stable} the frame $\R$ and the set of manifolds 
         $\Lambda=\{\Lambda_{j,\alpha}\,|\, (j,\alpha)\in J\}$ is a stable configuration near a point $\bar x\in \wtO$;

        \item[(5)']  for each double index $(j, \alpha)\in J$, $h_{j,\alpha}:
        \Lambda_{j,\alpha} \to \RR$ is bounded and $C^1$-smooth, with $h(\ox) \in \RR^{|J|}$;
        
         \item[(6)']  for each double index $(j,\alpha)\in J$ and each $k\in
        J_{j,\alpha}$, the right-hand side of \eqref{eq:final-system} is a
        function $\wtO \times \RR^{|J|} \to \RR$ which is bounded, continuous, and
        also uniformly Lipschitz in its second argument. 
 \end{itemize}
The statements (1)', (2)', (3)', and (5)' are self-evident, while statement (6)' 
follows from assumption (6) of Lemma \ref{lemma} together with 
assumptions (7) and (11) of the present theorem.

Finally, to verify (4)'  we show that stability of the configuration
$(\R,\Lambda)$ follows from that of the configuration $(\R,\Xi)$ (which holds
according to assumption (6) of Lemma \ref{lemma}).  
By stability of the configuration $(\R,\Xi)$, let $\widehat\Omega$  be a neighborhood of $\bar x$
for which there is,
for each $\alpha=1,\dots,m$, a {$\pa$-accessible} neighborhood 
$\Theta_{\alpha}\subset \RR^n$ such that the map 
$\psi_\alpha\colon \Theta_{\alpha} \to\widehat\Omega$ given by \eqref{diffeo} 
is a $C^1$-diffeomorphism (see Definitions~\ref{def-access}-\ref{def-stable}).
Then, for $(j,\alpha)\in J$,  the map
\begin{align*}
	\widehat\psi_{j,\alpha}(t_1,\dots,t_n)&:=
	W_{i_{\pa-j}}^{t_{\pa-j}} \cdots  W_{i_{1}}^{t_{1}}\, 
	\lambda_{j,\alpha}(t_{\pa-j+1},\dots,t_n)\\
	&\;=\psi_{\alpha}(t_{\pa-j+1},\dots,t_{\pa}, t_1,\dots,t_{\pa-j},t_{\pa+1},\dots,t_n )
\end{align*}
is defined  for all $(t_1,\dots,t_n)\in\RR^n$ with the property that
\[(t_{\pa-j+1},\dots,t_{\pa}, t_1,\dots,t_{\pa-j},t_{\pa+1},\dots,t_n )\in \Theta_\alpha.\]
Thus, the domain $\widehat \Theta_{j,\alpha}$ of $\widehat\psi_{j,\alpha}$ 
is the image of $\Theta_\alpha$ under a permutation-of-coordinates 
map $\Pi$ that interchanges the block of $1$st-through-$j$th coordinates 
with the block of $(j+1)$th-through-$(\pa-1)$th coordinates. The map  
$\widehat\psi_{j,\alpha}=\psi_\alpha\circ \Pi^{-1}$ is therefore a 
$C^1$-diffeomorphsim $\widehat \Theta_{j,\alpha}\to\widehat\Omega$. It is 
straightforward to verify that $(\pa-j)$-accessibility of $\widehat \Theta_{j,\alpha}$ 
follows from $\pa$-accessibility of $\Theta_\alpha$ (see Definition~\ref{def-access}.) 
This shows that assumption (4)' is satisfied.
   
Thus, according to Lemma~\ref{lemma},  there is a 
neighborhood $\wtO'$ of $\ox$ in $\RR^n$ on which there exists a unique 
set of functions $\{A_\ijj^\alpha(x) \,|\, (j,\alpha)\in J\}$ solving 
\eqref{eq:final-system} (with vanishing data on the $\Lambda_{j,\alpha}$),
for all $x\in \Lambda_{j,\alpha}^{k-j} \cap\wtO'$, 
where $\Lambda_{j,\alpha}^{k-j}$ are defined similarly to \eq{Xij} by
\beq\label{Lambdak}\Lambda_{j,\alpha}^{k-j}
	:=\{\widehat\psi_{j,\alpha}(t)\,|\, 
	t=(t_1,\dots,t_{k-j},0,\dots,0,t_{\pa-j+1},\dots,t_n)\in\widehat\Theta_{j,\alpha} \}
\eeq
Unwinding the  definition of $\widehat\psi_{j,\alpha}(t)$ we see that 
$\Lambda_{j,\alpha}^{k-j}$ is obtained by starting at a point $\Lambda_{j,\alpha}$, 
and then flowing in turn along the vectors fields whose indices appear no later than the 
$(k-j)$th member of $J_{j,\alpha}$. Recalling  that  $\Lambda_{j, \alpha}=\Xia^j$ and observing that  the $(k-j)$th member of $J_{j,\alpha}$ is the $k$-th member  of $I_\alpha$, we  conclude that $\Lambda_{j,\alpha}^{k-j}$ equals to $\Xia^k$, on which \eqref{eq:final-system} is to hold.

Finally,  since the identically vanishing functions $A_\ijj^\alpha(x)\equiv 0$ 
satisfy \eqref{eq:final-system} as well as the data, it follows from the uniqueness 
part of Lemma \ref{lemma} that the functions defined by \eq{A} are identically zero on  on $\wtO'$. That is,
    \[
        \br_\ijj(\tilde u_\alpha)\atx = \fa_\ijj(x,\wu(x))\qquad \text{for all $x\in \wtO'$,}
    \]
    as was to be shown.
    \end{proof}

\vskip5mm
\paragraph{\bf Acknowledgment:}This work was supported in part by the NSF grants  DMS-1311353 (PI: Jenssen) and DMS-1311743 (PI: Kogan).

\Addresses

\begin{bibdiv}
\begin{biblist}
\bib{mike_thesis}{thesis}{
   author={Benfield, Michael},
    title={Some Geometric Aspects of Hyperbolic Conservation Laws},
   date={2016},
   note={PhD thesis, NSCU, \url{https://repository.lib.ncsu.edu/handle/1840.16/11372}}
}
\bib{bjk_unpub}{article}{
   author={Benfield, Michael},
   author={Jenssen, Helge Kristian},
   author={Kogan, Irina A.},   
   title={On two integrability theorems of Darboux},
   date={2017},
   note={Unpublished manuscript,
   available for download from   \url{https://arxiv.org/abs/1709.07473}}
}

\bib{dar}{book}{
   author={Darboux, Gaston},
   title={Le\c cons sur les syst\`emes orthogonaux et les coordonn\'ees
   curvilignes. Principes de g\'eom\'etrie analytique},
   language={French},
   series={Les Grands Classiques Gauthier-Villars. [Gauthier-Villars Great
   Classics]},
   note={The first title is a reprint of the second (1910) edition; the
   second title is a reprint of the 1917 original;
   Cours de G\'eom\'etrie de la Facult\'e des Sciences. [Course on Geometry
   of the Faculty of Science]},
   publisher={\'Editions Jacques Gabay, Sceaux},
   date={1993},
   pages={600},
   isbn={2-87647-016-0},
   review={\MR{1365963}},
}

\bib{hart}{book}{
   author={Hartman, Philip},
   title={Ordinary differential equations},
   series={Classics in Applied Mathematics},
   volume={38},
   note={Corrected reprint of the second (1982) edition [Birkh\"auser,
   Boston, MA;  MR0658490 (83e:34002)];
   With a foreword by Peter Bates},
   publisher={Society for Industrial and Applied Mathematics (SIAM),
   Philadelphia, PA},
   date={2002},
   pages={xx+612},
   isbn={0-89871-510-5},
   review={\MR{1929104}},
   doi={10.1137/1.9780898719222},
}

\bib{lieb}{book}{
   author={Lieberstein, H. Melvin},
   title={Theory of partial differential equations},
   note={Mathematics in Science and Engineering, Vol. 93},
   publisher={Academic Press, New York-London},
   date={1972},
   pages={xiv+283},
   review={\MR{0355280}},
}
\bib{khar}{article}{
   author={Kharibegashvili, S.},
   title={Goursat and Darboux type problems for linear hyperbolic partial
   differential equations and systems},
   language={English, with English and Georgian summaries},
   journal={Mem. Differential Equations Math. Phys.},
   volume={4},
   date={1995},
   pages={127},
   issn={1512-0015},
   review={\MR{1415805}},
}
\bib{spi}{book}{
   author={Spivak, Michael},
   title={A comprehensive introduction to differential geometry. Vol. I},
   edition={2},
   publisher={Publish or Perish Inc.},
   place={Wilmington, Del.},
   date={1979},
   pages={xiv+668},
   isbn={0-914098-83-7},
   review={\MR{532830 (82g:53003a)}},
}
\bib{lang}{book}{
   author={Lang, Serge},
   title={Introduction to differentiable manifolds},
   series={Universitext},
   edition={2},
   publisher={Springer-Verlag, New York},
   date={2002},
   pages={xii+250},
   isbn={0-387-95477-5},
   review={\MR{1931083}},
}

\end{biblist}
\end{bibdiv}
\end{document}